\RequirePackage[l2tabu, orthodox]{nag}
\documentclass[10pt, draft]{article}

\usepackage{amsmath,amssymb,amsbsy,amsfonts,amsthm,latexsym,amsopn,amstext,
                            amsxtra,euscript,amscd}
\usepackage{booktabs}
\usepackage{enumitem}
\usepackage{tikz}
\usepackage[norefs,nocites]{refcheck}
\newtheorem{thm}{Theorem}
\newtheorem{lem}{Lemma}
\newtheorem{lemma}[thm]{Lemma}

\newtheorem{conj}[thm]{Conjecture}

\def\\{\cr}
\def\({\left(}
\def\){\right)}
\def\[{\left[}
\def\]{\right]}

\def\<{\langle}
\def\>{\rangle}
\def\fl#1{\left\lfloor#1\right\rfloor}

\def\cC{\mathcal C}
\def\cF{\mathcal F}

\def\cG{\mathcal G}

\def\notdivides{\mathrel{\kern-3pt\not\!\kern3.5pt\bigm|}}

\begin{document}


\title{\Large \textbf{Primes in floor function sets}}
\author{
\scshape {RANDELL HEYMAN}  \\
School of Mathematics and Statistics,\\ University of New South Wales \\
Sydney, Australia\\
\texttt {randell@unsw.edu.au}
}

\maketitle


\begin{abstract}
Let $x$ be a positive integer. We give an asymptotic formula for the number of primes in the set $\{\fl{x/n}, 1 \le n \le x\}$ and give some related results.
\end{abstract}


\section{Introduction}
There is an extensive body of research on
arithmetic functions with integer parts of real-valued functions, most commonly, with
 Beatty $\fl{\alpha n + \beta}$  sequences, see, for
 example,~\cite{ABS,BaBa,BaLi,GuNe,Harm}, and
 Piatetski--Shapiro $\fl{n^\gamma}$ sequences,
 see, for example,~\cite{Akb,BBBSW,BBGY,BGS,LSZ,Morg},
 with real $\alpha$, $\beta$ and $\gamma$.

 Recently there has been much research on sums of the form
 \begin{align}
 \label{eq:f sum}
 \sum_{n \le x}f\(\fl{\frac{x}{n}}\),
 \end{align}
 where throughout $x$ is a positive integer, $f$ is an arithmetic function and $\fl{\cdot}$ is the floor function.
 In \cite{Bor} the authors used exponential sums to find asymptotic bounds and formulas for various classes of arithmetic functions. Subsequent papers by various authors have mainly been focussed on improvements in exponential sums techniques (see \cite{Bor2,Che, Liu,Ma,Ma2,Stu,Wu,Zha,Zha2,Zhao}).

It is natural to examine more fundamental questions about the set$\{\fl{x/n}: 1 \le n \le x\}$. In \cite{Hey} an exact formula  for the cardinality of this set was given.
In this paper we count primes in this floor function set. Let
 $$\cG(x)=\left\{\fl{\frac{x}{n}}: 1 \le n \le x, \fl{\frac{x}{n}}\text{ is prime}\right\},$$
 and in particular $G(x):=|\cG(x)|$ This can estimated using exponential sums as follows:
\begin{thm}
 \label{thm:G(x)}
 We have
 $$G(x)=\frac{4 \sqrt{x}}{\log x}+O\(\frac{ \sqrt{x}}{(\log x)^2}\).$$
 \end{thm}

 The OEIS sequence A068050 attributes to Adams-Watters the statement that for $p$ prime not equal to 3 we have $G(p)=G(p-1)+1$. A proof does not seem evident. We prove the following:
\begin{thm}
\label{thm:adams}
Let $x$ be any prime not equal to 3. Then $G(x)=G(x-1)+1$.
\end{thm}
It is possible to link up $G(x)$ and $G(x-1)$ for some other classes of $x$ as follows:
 \begin{thm}
 \label{thm:G(x) semiprime}
Let $x=pq$ with $p,q$ odd primes, not necessarily distinct.
 Then
 $$G(x)=G(x-1)+1.$$
 \end{thm}
 These relationships between $G(x)$ and $G(x-1)$ may generalise, but with considerable difficulties.  For example, on a somewhat limited investigation using Maple we have the following:
 \begin{conj}
 \label{conj:G(x) 3 primes}
 Suppose $x=p_1p_2p_3$ with $2 <p_1< p_2 < p_3$.
 Then
 $$
G(x)= \begin{cases} G(x-1) & \text{if  } p_1p_2>p_3, \\
G(x-1)+1 & \text{if  } p_1p_2<p_3.
\end{cases}
$$
 \end{conj}

We can also examine the cardinality of the set
 $$\cF(x):=\left\{n:\fl{\frac{x}{n}}\text{ is prime}\right\}.$$
 This might more naturally be thought of as the cardinality of the subsequence $(\cF_{n_k})$ created from the sequence $(\cF_n)_{n=1}^x, \, \cF_n=\fl{x/n}$, where you retain $n$ for which $\cF_n$ is prime and remove $n$ for which $\cF_n$ is not prime. For example, we have $$\cF(10)=\{2,3,4,5\},$$ whilst it is more natural to think of the sequence (for $x=10$) $$(\cF_{n_k})=5,3,2,2.$$ Of course, the cardinalities are the same.
 The cardinality of $\cF(x)$ (or of $(\cF_{n_k})$) can be obtained by substituting $f(m)= \textbf{1}_\mathbb{P}(m)$ into \eqref{eq:f sum} and using recent results from Wu \cite{Wu} or from  Zhai \cite{Zha}. As is usual,  $\textbf{1}_\mathbb{P}(m)=1$ if $m$ is prime and 0 otherwise.
We obtain the following:
 \begin{thm}
\label{thm:n prime}
Let $F(x):=|\cF(x)|$.Then
$$F(x)=\mathcal{P}x+O\(x^{1/2}\),$$
where
$$\mathcal{P}=\sum_{n=1}^\infty \frac{\textbf{1}_\mathbb{P}(n)}{n(n+1)}=\sum_p \frac1{p(p+1)}\cong 0.330230.$$
\end{thm}
We can use an alternate elementary approach, without exponential sums, to arrive at a result with a slightly better lower bound. Specifically
\begin{thm}
\label{thm:n prime elementary}
There exists calculable constants $A_1$ and $A_2$ such that for all $x$,
$$\mathcal{P}x -\frac{A_1\sqrt{x}}{\log x} \le F(x) \le \mathcal{P}x+A_2\sqrt{x}.$$
\end{thm}

The methodology of Theorem \ref{thm:n prime} can be utilised for all indicator functions since these functions are all bounded by 1. For example, we state, but do not prove, the following:
\begin{thm}
We have
$$\left\{n:\fl{\frac{x}{n}}\text{ is a prime power}\right\}=\mathcal{D}x +O\(x^{1/2}\),$$
where
$$\mathcal{D}=\sum_{n=p^k} \frac1{n(n+1)}\cong 0.41382.$$
\end{thm}

Throughout we use $p$, with or without subscript, to denote a prime number.
The notation $f(x) = O(g(x))$  or $f(x) \ll g(x)$ is
equivalent to the assertion that there exists a constant $c>0$ such that $|f(x)|\le c|g(x)|$ for all $x$.
As is normal, we denote by $\Lambda$ the von Mangoldt function.

\section{Proof of Theorem \ref{thm:G(x)}}
We have
$$G(x)=\left|\left\{p: p \le x, p=\fl{\frac{x}{n}}\text{ for some } 1 \le n \le x\right\} \right|.$$
If $\fl{\frac{x}{n}}=p$ then
$$\frac{x}{p+1} < n \le \frac{x}{p}$$ and such an $n$ will exist if
$\fl{\frac{x}{p}}-\fl{\frac{x}{p+1}}>0$. So
$$G(x)=\sum_{p \le x} \delta\(\fl{\frac{x}{p}} -\fl{\frac{x}{p+1}}>0\),$$
where $\delta=1$ if the statement is true and 0 otherwise.
Let
\begin{align}
\label{eq:G}
G(x)&=G_1(x)+G_2(x)+G_3(x)+G_4(x),
\end{align}
where
$$G_1= \sum _{p < b}, \,\, G_2=\sum_{b \le p \le \sqrt{x}},\,\, G_3=\sum_{\sqrt{x} < p \le x^{34/67}},\,\, G_4=\sum_{x^{34/67} < p \le x},$$
and
$$b=\frac{\sqrt{4x+1}-1}{2}=\sqrt{x}+O(1).$$

For $G_1(x)$ the condition is always satisfied, since for $p < b$ we have
$$\fl{\frac{x}{p}}-\fl{\frac{x}{p+1}}>\frac{x}{p}-\frac{x}{p+1}-1=\frac{x}{p(p+1)}-1>0.$$
So
\begin{align}
\label{eq:G1}
G_1(x)&=\sum_{p \le b}1=\pi(\sqrt{x})+O(1)=\frac{2\sqrt{x}}{\log x}+O\(\frac{\sqrt{x}}{(\log x)^2}\).
\end{align}

Trivially
\begin{align}
\label{eq:G2}
G_2(x)&=O(1).
\end{align}

Next, we estimate $G_4(x)$. If $p > x^{34/67}$ then $p=\fl{\frac{x}{n}}$ for some $n \le x^{33/67}.$ Since there can be at most $x^{33/67}$ values for $n$ we have
\begin{align}
\label{eq:G4}
G_4(x)&=O\(x^{33/67}\).
\end{align}
For $G_3(x)$ (and $G_4(x)$) $p$ is large enough that $\fl{\frac{x}{p}}-\fl{\frac{x}{p+1}}$ can only equal 0 or 1.
So
$$G_3(x)=\sum_{\sqrt{x} < p \le x^{34/67}} \delta\(\fl{\frac{x}{p}} -\fl{\frac{x}{p+1}}>0\)=\sum_{\sqrt{x} < p \le x^{34/67}} \(\fl{\frac{x}{p}} -\fl{\frac{x}{p+1}}\).$$
Then, using $\psi(x)=x-\fl{x}-\frac1{2}$,
\begin{align}
\label{eq:G3}
G_3(x)&=x \sum_{\sqrt{x}<p\le x^{34/67}}\frac1{p(p+1)}+\sum_{\sqrt{x}<p\le x^{34/67}}\(\psi\(\frac{x}{p+1}\)-\psi\(\frac{x}{p}\)\).
\end{align}
Using partial summation and the Prime Number Theorem we have, for the first sum,
\begin{align}
\label{eq:G3 first}
x \sum_{\sqrt{x}<p\le x^{34/67}}\frac1{p(p+1)}&=x\sum_{\sqrt{x}<n\le x^{34/67}}\frac{\textbf{1}_\mathbb{P}(n)}{p(p+1)}\notag\\
&=
\end{align}

For the second sum of $G_3(x)$ we will use the following (\cite[Theorem 6.25]{Bor3}):
\begin{lem}
\label{lem:Bor}
Let $\delta \in [0,1], x \ge 1$ be a large real number and $R$, $R_1$ be positive integers such that $1 \le R \le R_1 \le 2R \le x^{2/3}$. Then , for all $\epsilon \in (0,\frac1{2}]$
$$x^{-\epsilon}\sum_{R \le n \le R_1} \Lambda(n) \psi\(\frac{x}{n+\delta}\)\ll \(x^2R^{33}\)^{1/38}+\(x^2R^{19}\)^{1/24}\(x^3R^2\)^{1/9}+\(x^3R^{-1}\)^{1/6}+R^{5/6}.$$
\end{lem}

Returning to the second sum of $G_3(x)$ and using the Lemma we have
\begin{align*}
\sum_{\sqrt{x}<p\le x^{34/67}}\(\psi\(\frac{x}{p+1}\)-\psi\(\frac{x}{p}\)\)
&\le \left|\sum_{\sqrt{x}<p\le x^{34/67}}\psi\(\frac{x}{p+1}\)\right|+ \left|\sum_{\sqrt{x}<p\le x^{34/67}}\psi\(\frac{x}{p}\)\right|.
\end{align*}
We now bound the sum involving $\psi(\frac{x}{p})$. The calculations for the sum involving $\psi(\frac{x}{p+1})$ is virtually identical. Let $m \sim N$ denote the inequalities $N < m \le 2N$.
We have
\begin{align*}
\sum_{\sqrt{x}<p\le x^{34/67}}\psi\(\frac{x}{p}\)&\ll \max_{\sqrt{x}<N \le x^{34/67}}\left|\sum_{p \sim N} \psi\(\frac{x}{p}\)\right|\log x.
\end{align*}
Next, using Abel summation,
\begin{align*}
\left|\sum_{p \sim N} \psi\(\frac{x}{p}\)\right|&=\left|\sum_{p \sim N}\(\frac1{\log p}\times \psi\(\frac{x}{p}\) \log p \)\right|\\
&\le \frac{2}{\log N} \max_{N \le N_2 \le N_1}\left|\sum _{N <p \le N_2} \psi\(\frac{x}{p}\)\log p+\textbf{1}_p(N) \psi\(
\frac{x}{N}\log N\)\right|\\
&\le \frac{2}{\log N} \max_{N \le N_2 \le N_1}\left\{\left|\sum_{N<n \le N_2}\Lambda(n)\psi\(\frac{x}{n}\)\right|+|R(N)|+\log N\right\},
\end{align*}
where
$$|R(N)| \le \sum_{\sqrt{N} < p \le \sqrt{N_2}}\log p \sum_{2 \le a \le \frac{\log N_2}{\log p}} 1<2 \sqrt{N}.$$
Using  Lemma \ref{lem:Bor} with $N_1=N_2=x^{34/67}$ and $N=\sqrt{x}$ we obtain
\begin{align}
\label{eq:G3 second}
\sum_{\sqrt{x}<p\le x^{34/67}}\psi\(\frac{x}{p}\)&\ll x^{\frac{1256}{2546}+\epsilon}.
\end{align}
Substituting \eqref{eq:G3 first} and \eqref{eq:G3 second} into \eqref{eq:G3} and we see that
$$G_3(x)=\frac{2 \sqrt{x}}{\log x}+O\(\frac{\sqrt{x}}{(\log x)^2}\),$$
and substituting this equation and
 \eqref{eq:G1}, \eqref{eq:G2} and \eqref{eq:G4} into \eqref{eq:G} completes the proof.

\section{Proof of Theorem \ref{thm:adams}}
Let $x$ be any prime not equal to 3. Suppose $p \in \cG(x)$ but $p\ne x$. So for some $n$ we have $\fl{x/n}=p$. As $x$ is a prime we have $x=np+u$ where $1 \le u \le n-1$. So $x-1=np+u-1$ from which
$$\fl{\frac{x-1}{n}}=\fl{\frac{np+u-1}{n}}=\fl{p+\frac{u-1}{n}}=p.$$
Thus $p \in \cG(x-1)$.

Conversely, suppose $p \in \cG(x-1)$.  So $\fl{(x-1)/n}=p$ for some $n$. So $x-1=np+u$ where $0 \le u \le n-1$. But $u \ne n-1$, for then $x=np+n=n(p+1)$, which contradicts the supposition that $x$ is prime. Thus $x-1=np+u$ with $0 \le u \le n-2$, and then
$$\fl{\frac{x}{n}}=\fl{\frac{np+u+1}{n}}=\fl{p+\frac{u+1}{n}}=p.$$
Therefore $p \in \cG(x)$.

We conclude that there is one-to-one correspondence between an $p\ne x \in \cG(x)$ and $p \ne x \in \cG(x-1)$.
Noting that we have $x \in \cG(x)$ but $x \not \in \cG(x-1)$ concludes the proof.

\section{Proof of Theorem \ref{thm:G(x) semiprime}}
We have $x=pq$ with  $p,q$ odd primes, not necessarily distinct. Without loss of generality assume $p \le q$.

\textit{Case 1:}Suppose that $r \in \cG(x)$ with $r \ne p,q$.
So $x=nr+u$ with $0 \le u \le n-1$. But if $u=0$ then $x=nr$ which is impossible since $r \ne p,q$. So $1 \le u \le n-1$ and thus
$$\fl{\frac{x-1}{n}}=\fl{\frac{nr+u-1}{n}}=\fl{r+\frac{u-1}{n}}=r.$$
So $r \in \cG(x)$.
Conversely, suppose $r \in \cG(x-1)$. So $x-1=nr+u$ with $0 \le u \le n-1$. But if $u=n-1$ then $x=nr+n-1+1=(n+1)r$ which is again impossible since $r \ne p,q$. So $0 \le u \le n-2$ and then
$$\fl{\frac{x}{n}}=\fl{\frac{nr+u+1}{n}}=\fl{r+\frac{u+1}{n}}=r.$$
So $r \in \cG(x-1)$.

\textit{Case 2:} Suppose that $r \in \cG(x)$ with $r=p=q$.
Since $r=\fl{x/r}$ we have $r \in \cG(x)$. But
$$\fl{\frac{x-1}{r}}=\fl{r-\frac1{r}}=r-1$$
and
$$\fl{\frac{x-1}{r-1}}=\fl{\frac{r^2-1}{r-1}}=r+1.$$
So $r \not\in \cG(x-1).$

\textit{Case 3:} Suppose that $r \in \cG(x)$ with $r=p \ne q$ or $r=q \ne p$.
Recall $p \le  q$. It is clear that $p \in \cG(x)$ and $q \in \cG(x)$.
Then
$$\fl{\frac{x-1}{q-1}}=\fl{\frac{pq-1}{q-1}}=\fl{p+\frac{p-1}{q-1}}=p,$$
so $p \in \cG(x-1)$.
But
$$\fl{\frac{x-1}{p}}=\fl{\frac{pq-1}{p}}=\fl{q-\frac{1}{p}}=q-1,$$
and
$$\fl{\frac{x-1}{p-1}}=\fl{\frac{pq-1}{p-1}}=\fl{q+\frac{q-1}{p-1}}>q.$$
So $q \not \in \cG(x-1)$.

Reviewing the three cases we see that $G(x)=G(x-1)+1$, which proves the theorem.

\section{Proof of Theorem \ref{thm:n prime}}
We have
$$F(x)=\sum_{n \le x} \textbf{1}_\mathbb{P}\(\fl{\frac{x}{n}}\).$$
We will require the following, proven independently by Wu \cite[Theorem]{Wu} and Zhai \cite[Theorem 1]{Zha}:
\begin{lemma}
\label{thm:Sfx}
Let $f$ be a complex-valued arithmetic function such that $f(n) \ll n^{\alpha}(\log n)^{\theta}$ for some $\alpha \in [0,1)$ and $\theta \ge 0$.  Then
$$\sum_{n \leqslant x} f \(\fl{x/n} \) = x \sum_{n=1}^{\infty} \frac{f(n)}{n(n+1)} + O \( x^{\frac{1}{2}(\alpha+1)}(\log x)^{\theta}\).$$
\end{lemma}
Using this lemma with $\alpha=0$ and  $\theta=0$ we have
$$F(x)=\mathcal{P}x+O\(x^{1/2}\),$$
where
$$\mathcal{P}=\sum_{n=1}^\infty \frac{\textbf{1}_\mathbb{P}(n)}{n(n+1)}=\sum_p \frac1{p(p+1)}\cong 0.330230,$$
completing the proof.

We note, in passing, that
$$\sum_p \frac1{p(p+1)}=\sum_{s=2}^\infty \sum_p \frac{(-1)^s}{p^s}=\sum_{s=2}^\infty\sum_{n=1}^\infty \mu(n) \frac{\log \zeta(ns)}{n}.$$

\section{Proof of Theorem \ref{thm:n prime elementary}}
Let $\mathbb{P}$ be the set of (positive) primes and $\overline{\mathbb{P}}$ be the set of (positive) non-primes. We create upper and lower bounds for $F(x)$ from the set
$$\cC(x):=\{n: 1 \le n \le x\}.$$
The set $\cF(x)$ is simply the set $\cC(x)$ after removing all $n$ such that $\fl{x/n}$ is a non-prime from $\cC(x)$.
For the upper bound we truncate the process by removing from $\cC(x)$ only those $n$ such that $\fl{x/n}$ is a non-prime less than or equal to $\sqrt{x}$.
In total we remove
$$\sum_{\substack{c \in \overline{\mathbb{P}}\\c \le \sqrt{x}}}\(\fl{\frac{x}{c}}-\fl{\frac{x}{c+1}}\)>\sum_{\substack{c \in \overline{\mathbb{P}}\\c \le \sqrt{x}}}\(\frac{x}{c(c+1)}-1\)$$
values of $n$. So
\begin{align*}
F(x)&< x-\sum_{\substack{c \in \overline{\mathbb{P}}\\c \le \sqrt{x}}}\(\frac{x}{c(c+1)}-1\)\\
&= x-\sum_{\substack{c \in \overline{\mathbb{P}}\\c \le \sqrt{x}}}\frac{x}{c(c+1)}+\sum_{\substack{c \in \overline{\mathbb{P}}\\c \le \sqrt{x}}}1\\
&=x-x\(\sum_{c \le \sqrt{x}} \frac1{c(c+1)}-\sum_{\substack{c \in \mathbb{P}\\c \le \sqrt{x}}}\frac{1}{c(c+1)}\)+\sqrt{x}-\pi(\sqrt{x})\\
&=x-\frac{x\sqrt{x}}{\sqrt{x}+1} +x\sum_{\substack{c \in \mathbb{P}\\c \le \sqrt{x}}}\frac{1}{c(c+1)}+O\(\sqrt{x}\)\\
&=\frac{x}{\sqrt{x}+1}+x\sum_{p}\frac1{p(p+1)} -x\sum_{p>\sqrt{x}}\frac1{p(p+1)} +O\(\sqrt{x}\).
\end{align*}
Then
$$\sum_{p>\sqrt{x}}\frac1{p(p+1)}\le \sum_{n>\sqrt{x}}\frac1{n(n+1)}=O\(\frac1{\sqrt{x}}\),$$
and so
\begin{align}
F(x)& \le \mathcal{P}x+O\(\sqrt{x}\).
\end{align}
For the lower bound we add up the number of $n$ for which $\fl{x/n}$ is a prime less than or equal to $\sqrt{x}$.
Then, using $\pi(m)=\frac{m}{\log m} + O\(\frac{m}{(\log m)^2}\)$ and Riemann-Stieltjes integration,
\begin{align*}
F(x)&=\sum_{\substack{c \in \mathbb{P}\\c \le \sqrt{x}}}\(\fl{\frac{x}{c}}-\fl{\frac{x}{c+1}}\)
\ge \sum_{\substack{c \in \mathbb{P}\\c \le \sqrt{x}}} \(\frac{x}{c(c+1)}+1\)\\
&=x\sum_{\substack{c \in \mathbb{P}\\c \le \sqrt{x}}} \frac{1}{c(c+1)}+\sum_{\substack{c \in \mathbb{P}\\c \le \sqrt{x}}}1
=x\(\mathcal{P}-\sum_{p > \sqrt{x}}\frac{1}{p(p+1)}\)+\sum_{\substack{c \in \mathbb{P}\\c \le \sqrt{x}}}1\\
&=\mathcal{P}x-x \sum_{p>\sqrt{x}}\frac1{p(p+1)}+O\(\frac{\sqrt{x}}{\log x}\)=\mathcal{P}x+O\(\frac{\sqrt{x}}{\log x}\),
\end{align*}
completing the proof.

\section{Acknowledgements}
I thank Joshua Stucky who sketched out the proof of the asymptotic formula for $G(x)$. I thank Olivier Bordell\`es for assisting in the proof of Theorem \ref{thm:G(x)}. I thank Igor Shparlinski for some initial discussions. I thank Stan Wagon for pointing out that some of the inequalities in Conjecture 4 should be changed.

\makeatletter
\renewcommand{\@biblabel}[1]{[#1]\hfill}
\makeatother

\end{document}